\documentclass[a4paper,12pt]{article}

\usepackage[english]{babel} 
\usepackage[ansinew]{inputenc}
\usepackage[LGR,T1]{fontenc}

\usepackage{verbatim}
\usepackage{color}
\usepackage{calc}
\usepackage{hyperref}
\usepackage[dvips]{graphicx}
\hypersetup{colorlinks=true, linkcolor=blue, filecolor=blue, urlcolor=blue}

\usepackage{amssymb}
\usepackage[mathscr]{eucal}
\usepackage{amscd,amsmath,amsfonts,amssymb,textcomp}
\everymath{\displaystyle}

\usepackage[all]{xy}

\newcommand{\MBC}{\mathbb{C}}

\newcommand{\MBN}{\mathbb{N}}

\newcommand{\MBQ}{\mathbb{Q}}

\newcommand{\MBZ}{\mathbb{Z}}

\newcommand{\MCC}{\mathcal{C}}

\newcommand{\MCE}{\mathcal{E}}

\newcommand{\MCG}{\mathcal{G}}
\newcommand{\MCH}{\mathcal{H}}
\newcommand{\MCI}{\mathcal{I}}

\newcommand{\MCL}{\mathcal{L}}

\newcommand{\MCO}{\mathcal{O}}

\newcommand{\MCQ}{\mathcal{Q}}

\newcommand{\MCS}{\mathcal{S}}
\newcommand{\MCT}{\mathcal{T}}
\newcommand{\MCU}{\mathcal{U}}

\newcommand{\MFa}{\mathfrak{a}}

\newcommand{\MFc}{\mathfrak{c}}

\newcommand{\MFf}{\mathfrak{f}}

\newcommand{\MFK}{\mathfrak{K}}

\newcommand{\MFm}{\mathfrak{m}}

\newcommand{\MFp}{\mathfrak{p}}
\newcommand{\MFP}{\mathfrak{P}}
\newcommand{\MFq}{\mathfrak{q}}
\newcommand{\MFQ}{\mathfrak{Q}}
\newcommand{\MFr}{\mathfrak{r}}

\newcommand{\MFu}{\mathfrak{u}}

\newcommand{\MSU}{\mathscr{U}}

\newcommand{\SCm}{\textsc{m}}

\newcommand{\SCr}{\textsc{r}}

\newcommand{\GGa}{\alpha}
\newcommand{\GGb}{\beta}

\newcommand{\GGD}{\Delta}
\newcommand{\GGe}{\epsilon}
\newcommand{\GVe}{\varepsilon}

\newcommand{\GVf}{\varphi}
\newcommand{\GGg}{\gamma}
\newcommand{\GGG}{\Gamma}
\newcommand{\GGi}{\iota}
\newcommand{\GGk}{\kappa}

\newcommand{\GGl}{\lambda}
\newcommand{\GGL}{\Lambda}
\newcommand{\GVp}{\varpi}

\newcommand{\GGs}{\sigma}

\newcommand{\GGt}{\theta}
\newcommand{\GGT}{\Theta}
\newcommand{\GVt}{\vartheta}
\newcommand{\GGu}{\upsilon}

\newcommand{\GGw}{\omega}
\newcommand{\GGW}{\Omega}
\newcommand{\GGz}{\zeta}

\newcommand{\ba}[2]{{\mathop{#1}\limits_{#2}}}
\newcommand{\op}[3]{{\mathop{#1}\limits_{#2}^{#3}}}

\newcommand{\lA}{\left\{}
\newcommand{\rA}{\right\}}
\newcommand{\la}{\left\langle}
\newcommand{\ra}{\right\rangle}

\newcommand{\Ann}{\mathrm{Ann}}
\newcommand{\Char}{\mathrm{char}}
\newcommand{\cl}{\mathrm{cl}}
\newcommand{\Cl}{\mathrm{Cl}}

\newcommand{\Cok}{\mathrm{Cok}}

\newcommand{\Gal}{\mathrm{Gal}}
\newcommand{\Hom}{\mathrm{Hom}}

\newcommand{\IM}{\mathrm{Im}}
\newcommand{\Ker}{\mathrm{Ker}}

\newcommand{\PR}{\mathrm{pr}}

\newcommand{\Tor}{\mathrm{Tor}}
\newcommand{\Tr}{\mathrm{Tr}}

\usepackage{makeidx}
\title{ON THE CLASSICAL MAIN CONJECTURE FOR IMAGINARY QUADRATIC FIELDS}
\author{\textsc{St\'ephane VIGUI\'E}
\footnote{St\'ephane Vigui\'e, Laboratoire de math\'ematiques de Besançon, UMR CNRS 6623, Universit\'e de Franche-Comt\'e, 16 route de Gray, 25030 Besançon cedex, France.
e-mail: \texttt{stephane.viguie@univ-fcomte.fr}}
}
\makeindex

\newtheorem{dfe}{dfe}[section]
\newtheorem{lem}{lem}[section]
\newtheorem{pro}{pro}[section]
\newtheorem{rem}{rem}[section]
\newtheorem{teh}{teh}[section]

\newtheorem{df}[dfe]{Definition}
\newtheorem{lm}[lem]{Lemma}
\newtheorem{pr}[pro]{Proposition}

\newtheorem{rmq}[rem]{Remark}
\newtheorem{theo}[teh]{Theorem}

\numberwithin{equation}{section}

\begin{document}
\maketitle

\begin{abstract}
Let $p$ be a prime number, and let $k$ be an imaginary quadratic number field in which $p$ decomposes into two distinct primes $\MFp$ and $\bar{\MFp}$.
Let $k_\infty$ be the unique $\MBZ_p$-extension of $k$ which is unramified outside of $\MFp$, and let $K_\infty$ be a finite extension of $k_\infty$, abelian over $k$.
In case $p\notin\{2,3\}$, we prove that in $K_\infty$, the characteristic ideal of the projective limit of the $p$-class group coincides with the characteristic ideal of the projective limit of units modulo elliptic units. 
Our approach is based on Euler systems, which were first used in this context by Rubin in \cite{rubin88}.
For $p\in\{2,3\}$, we obtain a divisibility relation, up to a certain constant.
\end{abstract}
\bigskip

\\


\section{Introduction.}

Let $p$ be a prime number, and let $k$ be an imaginary quadratic number field in which $p$ decomposes into two distinct primes $\MFp$ and $\bar{\MFp}$.
Let $k_\infty$ be the unique $\MBZ_p$-extension of $k$ which is unramified outside of $\MFp$, and let $K_\infty$ be a finite extension of $k_\infty$, abelian over $k$.
Let $G_\infty$ be the Galois group of $K_\infty/k$.
We choose a decomposition of $G_\infty$ as a direct product of a finite group $G$ (the torsion subgroup of $G_\infty$) and a topological group $\GGG$ isomorphic to $\MBZ_p$, $G_\infty = G\times \GGG$.
For any $n\in\MBN$, let $K_n$ be the field fixed by $\GGG_n := \GGG^{p^n}$, and let $G_n:=\Gal\left(K_n/k\right)$.
Remark that there may be different choices for $\GGG$, but when $p^n$ is larger than the order of the $p$-part of $G$, the group $\GGG_n$ does not depend on the choice of $\GGG$.

Let $F/k$ be an abelian extension of $k$.
If $[F:k]<\infty$, we denote by $\MCO_F$ the ring of integers of $F$.
We write $\MCO_F^\times$ for the group of global units of $F$, and $C_F$ for the group of elliptic units of $F$ (see section \ref{sectionelliptic}).
Also we let $A_F$ be the $p$-part of the class group $\Cl\left(\MCO_F\right)$ of $\MCO_F$.
We set $\MCE_F:=\MBZ_p\otimes_\MBZ\MCO_F^\times$ and $\MCC_F:=\MBZ_p\otimes_\MBZ C_F$.
When $F/k$ is infinite, we define $\MCE_F$, $\MCC_F$ and $A_F$, by taking projective limits over finite sub-extensions, under the norm maps.
For any $n\in\MBN\cup\{\infty\}$, we set $\MCE_n:=\MCE_{K_n}$, $\MCC_n:=\MCC_{K_n}$, and $A_n:=A_{K_n}$.

For any profinite group $\MCG$, and any commutative ring $R$, we define the Iwasawa algebra
\[R\left[\left[\MCG\right]\right] := \varprojlim R\left[\MCH\right],\]
where the projective limit is over all finite quotient $\MCH$ of $\MCG$.
In case $\MCG=G_\infty$, we shall write
\[\GGL := \MBZ_p\left[\left[G_\infty\right]\right].\]
Then $A_\infty$ and $\MCE_\infty/\MCC_\infty$ are naturally $\GGL$-modules.
As we shall see below, they are finitely generated and torsion over $\MBZ_p[[\GGG]]$.
Let us fix a topological generator $\GGg$ of $\GGG$, and set $T:=\GGg-1$.
Then for any finite extension $L/\MBQ_p$, $\MCO_L[[\GGG]]$ is isomorphic to $\MCO_L[[T]]$, where $\MCO_L$ is the ring of integers of $L$.
It is well known that $\MCO_L[[T]]$ is a noetherian, regular, local domain.
We also recall that $\MCO_L[[T]]$ is a unique factorization domain, and that $\MBQ_p\otimes_{\MBZ_p}\MCO_L[[T]]$ is a principal ring.
If $u$ is a uniformizer of $\MCO_L$, then the maximal ideal $\MFm$ of $\MCO_L$ is generated by $u$ and $T$, and $\MCO_L[[T]]$ is a compact topological ring with respect to its $\MFm$-adic topology.
A morphism $f:M\rightarrow N$ between two finitely generated $\MCO_L[[T]]$-module is called a pseudo-isomorphism if its kernel and its cokernel are finite.
If a finitely generated $\MCO_L[[T]]$-module $M$ is given, then one may find elements $P_1$, ..., $P_r$ in $\MCO_L[T]$, irreducible in $\MCO_L[[T]]$, and nonnegative integers $n_0$, ..., $n_r$, such that there is a pseudo-isomorphism
\[M\longrightarrow \MCO_L[[T]]^{n_0} \oplus \bigoplus_{i=1}^r \MCO_L[[T]]/\left(P_i^{n_i}\right).\]
Moreover, the integer $n_0$ and the ideals $\left(P_1^{n_1}\right)$, ..., $\left(P_r^{n_r}\right)$, are uniquely determined by $M$.
If $n_0=0$, then the ideal generated by $P_1^{n_1}\cdots P_r^{n_r}$ is called the characteristic ideal of $M$, and is denoted by $\Char_{\MCO_L[[T]]}(M)$.

We denote by $\MBC_p$ a completion of an algebraic closure of $\MBQ_p$.
Let $\chi:G\rightarrow\MBC_p$ be an irreducible character of $G$.
Let $\MBQ_p(\chi)\subset\MBC_p$ be the abelian extension of $\MBQ_p$ generated by the values of $\chi$.
We denote by $\MBZ_p(\chi)$ the ring of integers of $\MBQ_p(\chi)$.
The group $G$ acts naturally on $\MBQ_p(\chi)$.
We recall that if $g\in G$ and $x\in\MBQ_p(\chi)$ then $g.x:=\chi(g)x$.
For any $\MBZ_p[G]$-module $Y$, we define its $\chi$-quotient $Y_\chi:=\MBZ_p(\chi)\otimes_{\MBZ_p[G]}Y$.
Moreover, if we set 
\[\psi:G\longrightarrow\MBQ_p,\quad\GGs \mapsto \Tr \left(\chi(\GGs)\right),\] 
where $\Tr$ is the trace map for the extension $\MBQ_p(\chi)/\MBQ_p$, then $\psi$ is an irreducible $\MBQ_p$ character on $G$, and we will write $\chi\vert\psi$.
Recall that any irreducible $\MBQ_p$ character on $G$ can be obtained in this way.
We define the idempotent $e_\psi$ of $\MBQ_p[G]$ attached to $\psi$,
\[e_\psi\,= \frac{1}{\#G}\,\sum_{g\in G}\,\psi(g)\,g^{-1}.\]
The restriction to the $\psi$-part of the canonical surjective map $\MBQ_p\otimes_{\MBZ_p}Y \rightarrow \MBQ_p\otimes_{\MBZ_p}Y_\chi$ is an isomorphism of $\MBQ_p[G]$-modules,
\begin{equation}
\label{isopsich}
e_\psi \left(\MBQ_p\otimes_{\MBZ_p}Y\right) \simeq \MBQ_p\otimes_{\MBZ_p}Y_\chi \simeq \MBQ_p(\chi)\otimes_{\MBZ_p[G]}Y.
\end{equation}
Also, we have the following decomposition
\[\MBQ_p\otimes_{\MBZ_p}Y\; = \; \mathop{\oplus}_\psi\, e_\psi\,\left(\MBQ_p\otimes_{\MBZ_p}Y\right),\]
where the sum is over all irreducible $\MBQ_p$ characters $\psi$ on $G$.
Finally, if $Y$ is a $\GGL$-module, then $Y_\chi$ is a $\MBZ_p(\chi)[[T]]$-module in a natural way.
As a particular case, $\GGL_\chi \simeq \MBZ_p(\chi)[[T]]$.
For any finitely generated $\GGL_\chi$-module $Z$, we shall denote $\Char_{\GGL_\chi}Z$ simply by $\Char Z$.

The goal of this article is to prove Theorem \ref{mainconj} below, which is a formulation of the (one-variable) main conjecture.
In \cite[Theorem 4.1]{rubin91} and \cite[Theorem 2]{rubin94}, Rubin used Euler systems to prove the main conjectures for $\MBZ_p$ or $\MBZ_p^2$ extensions of a finite abelian extension $F$ of $k$, where $p\nmid w_k[F:k]$, $w_k$ being the number of roots of unity in $k$.
More recently, Hassan Oukhaba adapted Rubin's method and obtained Theorem \ref{mainconj} for $p=2$, still under the condition $2\nmid\left[K_0:k\right]$ (see \cite{oukhaba10}).
Inspired by the ideas of Rubin, Greither used Euler systems to prove the main conjecture for cyclotomic units and for the cyclotomic $\MBZ_p$-extension $F_\infty/F$, with $F_\infty$ abelian over $\MBQ$ (see \cite[Theorem 3.2]{greither92}).
Bley proved Theorem \ref{mainconj} when $p\nmid2\#\left(\Cl\left(\MCO_k\right)\right)$, and when there is a nonzero ideal $\MFf$ of $\MCO_k$, prime to $\MFp$, 
such that for all $n\in\MBN$, $K_n = k\left(\MFf\MFp^{n+1}\right)$ is the ray class field of $k$ modulo $\MFf\MFp^{n+1}$ (see \cite[Theorem 3.1]{bley06}).
Here we prove the general case.

Also, we draw the attention of the reader to a cohomological two-variables main conjecture, which has been recently proved for all primes by J.\,Johnson-Leung and G.\,Kings in \cite{johnson-leung-kings}, as a consequence of the Tamagawa number conjecture.
In their treatment they replaced $\chi$-quotients by Galois cohomology with coefficients in the Galois representations defined by $\chi$, and used Euler systems as defined by Kato.
From their result, they deduced the classical two-variables main conjecture for $\MBZ_p^2$-extensions $F_\infty := \mathop{\cup}_{n=0}^\infty k\left(p^n\MFf\right)$ 
where $\MFf$ is any nonzero ideal of $\MCO_k$, and when $p$ does not divide the torsion subgroup of $\Gal\left(F_\infty/k\right)$.

\begin{theo}\label{mainconj}
Let $\chi$ be an irreducible $\MBC_p$ character on $G$.

$\mathrm{(i)}$ If $p\notin\{2,3\}$, then $\Char\left(A_{\infty,\chi}\right) = \Char\left(\MCE_\infty/\MCC_\infty\right)_\chi$.

$\mathrm{(ii)}$ If $p\in\{2,3\}$, then there is $m_\chi\in\MBN$ such that 
\begin{equation}
\label{divisibility}
\Char\left(A_{\infty,\chi}\right) \quad\text{divides}\quad \MFu_\chi^{m_\chi} \Char\left(\MCE_\infty/\MCC_\infty\right)_\chi,
\end{equation}
where $\MFu_\chi$ is a uniformizer of $\MBZ_p(\chi)$.
\end{theo}

\section{Semi-local units.}\label{sectionellipticunits}

For every $n\in\MBN$, we denote by $\MCU_n$ the $\MBZ_p[G_n]$-module of principal semi-local units over the primes above $\MFp$.
We define 
\[\MCU_\infty \; := \; \varprojlim\,\MCU_n,\]
by taking the projective limit under the norm maps.

For any $n\in\MBN$, we write $\GGg_n$ for $\GGg^{p^n}$.
Then for any $\MBZ_p[[T]]$-module $M$, we denote by $M^{\GGG_n}$ the module of $\GGG_n$-invariants of $M$, and we denote by $M_{\GGG_n}$ the module of $\GGG_n$-coinvariants of $M$.
By definition, they are respectively the kernel and the cokernel of the multiplication by $1-\GGg_n$ on $M$.

\begin{pr}\label{QptensorZpsemilocalunits}
$\MBQ_p \otimes_{\MBZ_p} \MCU_\infty$ is a free $\MBQ_p \otimes_{\MBZ_p} \GGL$-module of dimension $1$.
\end{pr}

\noindent\textsl{Proof.}
Let $\MFP$ be a prime of $K_\infty$ above $\MFp$, and let $k'$ be the completion of $k$ at $\MFp$.
For every $n\in\MBN$ we respectively denote by $K'_n$, $\MCO'_n$, and $\widehat{\MFP}_n$ the completion of $K_n$ at $\MFP$, 
the ring of integers of $K'_n$, and the maximal ideal of $\MCO'_n$.
Let us also denote the group $1 + \widehat{\MFP}_n$ by $\MCU'_n$.
For sufficiently large $m$, the $p$-adic logarithm is an isomorphism of $\MBZ_p\left[\Gal\left(K'_n / k'\right)\right]$-modules from $1+\widehat{\MFP}_n^m$ into $\widehat{\MFP}_n^m$.
Taking the tensor product of this isomorphism by $\MBQ_p$ over $\MBZ_p$, we see that
\begin{equation}
\label{tralalalal}
\MBQ_p \otimes_{\MBZ_p} \MCU'_n \; \simeq \; K'_n \; \simeq \; \MBQ_p\left[ \Gal\left(K'_n / k'\right) \right],
\end{equation}
where the last isomorphism holds by the normal basis theorem, and since $k'=\MBQ_p$.
The field $K'_\infty := \mathop{\cup}_{n=0}^\infty K'_n$ is a $\MBZ_p$-extension of $K'_0$, abelian over $k'$.
The Galois group of $K'_\infty / k'$ is canonically identified to the decomposition group of $\MFp$ in $K_\infty / k$.
Hence we only have to show that the projective limit
\[\MCU'_{\infty} \; := \; \varprojlim \, \MCU'_n,\]
with respect to the norm maps, is such that $\MBQ_p \otimes_{\MBZ_p} \MCU'_{\infty}$ is a free $\MBQ_p \otimes_{\MBZ_p} \GGL'$-module, where 
\[\GGL' \; := \; \MBZ_p \left[\left[ \Gal\left(K'_{\infty} / k'\right) \right]\right].\]
We have a decomposition $\Gal\left(K'_{\infty} / k'\right) = G'\times \GGG'$, where $G'\subseteq G$, and $\GGG'\simeq\MBZ_p$ as a pro-$p$-group.
Remark that for $n$ large enough, we have $\GGG_n\subseteq\GGG'$.
Let $\chi$ be any $\MBC_p$ irreducible character on $G'$, and let $\psi : G' \rightarrow \MBQ_p$ be the unique irreducible character such that $\chi\vert\psi$.
We must prove that $e_\psi\left( \MBQ_p \otimes_{\MBZ_p} \MCU'_\infty \right)$ is a free $e_\psi\left( \MBQ_p \otimes_{\MBZ_p} \GGL' \right)$-module.
Since the group of $p$-power roots of unity in $K'_\infty$ is finite by Lemma \ref{mufini} below, we deduce from \cite[Theorem 25]{iwasawa73} 
that $\MCU'_\infty\hookrightarrow\MBZ_p\left[\left[\GGG'\right]\right]^d$ where $d=\left[K'_0:k'\right]$.
In particular $e_\psi \left(\MBQ_p\otimes_{\MBZ_p}\MCU'_\infty\right) \subseteq \left(\MBQ_p\otimes_{\MBZ_p}\MBZ_p\left[\left[\GGG'\right]\right]\right)^d$ 
and hence is a torsion free $e_\psi\left(\MBQ_p\otimes_{\MBZ_p}\GGL'\right)$-module.
Let us remark that $e_\psi\left(\MBQ_p\otimes_{\MBZ_p}\GGL'\right) \simeq \MBQ_p \otimes_{\MBZ_p} \MBZ_p(\chi)\left[\left[\GGG'\right]\right]$ is a principal ring.
We deduce that $e_\psi\MBQ_p\otimes_{\MBZ_p}\MCU'_{\infty}$ is a free $e_\psi \left( \MBQ_p \otimes_{\MBZ_p} \GGL \right)$-module. 
Let $r_\chi$ be its rank.
We must show that $r_\chi=1$.
Let us choose $n$ such that $\GGG_n\subseteq\GGG'$.
Then
\begin{equation}
e_\psi \left( \MBQ_p\otimes_{\MBZ_p} \left(\MCU'_{\infty}\right)_{\GGG_n} \right) \simeq e_\psi \left( \MBQ_p \otimes_{\MBZ_p} \MBZ_p \left[ \Gal\left( K'_n/k \right) \right] \right) ^{r_\chi}.
\label{kjgkgh}
\end{equation} 
Exactly as in \cite[Proof of Proposition 3.6]{oukhaba10}, one can prove that the kernel and the cokernel of the canonical map $\left(\MCU'_{\infty}\right)_{\GGG_n}\rightarrow\MCU'_n$ are finitely generated $\MBZ_p$-modules of rank $1$, and invariant under the action of $\Gal\left(K'_{\infty} / k'\right)$.
We deduce that
\begin{equation}
\label{JGHGHGGH}
\dim_{\MBQ_p} \left( e_\psi \left( \MBQ_p\otimes_{\MBZ_p} \left(\MCU'_{\infty}\right)_{\GGG_n} \right) \right) = \dim_{\MBQ_p} \left( e_\psi \left( \MBQ_p \otimes_{\MBZ_p} \MCU'_n \right) \right).
\end{equation}
But $e_\psi \left( \MBQ_p \otimes_{\MBZ_p} \MCU'_n \right) \simeq e_\psi \left( \MBQ_p \otimes_{\MBZ_p} \MBZ_p \left[ \Gal\left( K'_n/k' \right) \right] \right)$ by (\ref{tralalalal}).
Thus $r_\chi=1$ by (\ref{JGHGHGGH}) and (\ref{kjgkgh}).
\hfill $\square$
\\

\begin{lm}\label{mufini}
Let the notation be as in the proof of Proposition \ref{QptensorZpsemilocalunits}.
Then the group $\mu_{p^\infty} \left( K'_{\infty} \right)$ of $p$-power roots of unity in $K'_\infty$ is finite.
\end{lm}

\noindent\textsl{Proof.}
As previously, we write $k'$ for the completion of $k$ at $\MFp$.
Since $k'=\MBQ_p$, it is well known that the kernel of the local norm residue symbol 
\[\left(\cdot, k'\left( \mu_{p^\infty} \right) / k' \right) : \left(k'\right)^\times \rightarrow \Gal\left(k'\left( \mu_{p^\infty} \right) / k' \right)\] 
is the free group $\la p\ra$ generated by $p$ (see for instance \cite[p.\,323, Proposition (1.8)]{neukirch99}).
Assume $\mu_{p^\infty}\subset K'_{\infty}$.
Then the kernel of the local norm residue symbol 
\[\left(\cdot, K'_{\infty} / k' \right) : \left(k'\right)^\times \rightarrow \Gal\left(K'_{\infty} / k' \right)\]
is a subgroup of $\la p\ra$, whose index is finite.
Let $\MFQ$ be a prime of $k_\infty$ above $\bar{\MFp}$.
We write $k''$ for the completion of $k$ at $\bar{\MFp}$.
For all $n\in\MBN$, we denote by $k'_n$ (resp. $k''_n$) the completion of $k_n$ at $\MFP$ (resp. $\MFQ$).
We set $k'_\infty := \mathop{\cup}_{i=0}^nk'_n$ and $k''_\infty := \mathop{\cup}_{i=0}^nk''_n$.
Since $\bar{\MFp}$ is finitely decomposed in $k_\infty / k$, the extension $k''_{\infty} / k''$ is infinite.
But it is also unramified, and then its Galois group is topologically generated by $\left( p, k''_{\infty} / k'' \right)$.
By the product formula, and since $k_\infty/k$ is unramified outside of $\MFp$, we have 
$\left( p, k''_{\infty} / k'' \right)_{|k_\infty} = \left( p^{-1}, k'_{\infty} / k' \right)_{|k_\infty}$, 
and we deduce that for all $n\in\MBZ\setminus\{0\}$, $\left( p^n, K'_{\infty} / k' \right)\neq1$.
Hence $\left(\cdot, K'_{\infty} / k' \right)$ is injective, which is absurd.
\hfill $\square$
\\

\section{Elliptic units.}\label{sectionelliptic}

For $L$ and $L'$ two $\MBZ$-lattices of $\MBC$ such that $L\subseteq L'$ and $[L':L]$ is prime to $6$, we denote by $z\mapsto\psi\left(z;L,L'\right)$ the elliptic function defined in \cite{robert90}.
For $\MFm$ a nonzero proper ideal of $\MCO_k$, and $\MFa$ a nonzero ideal of $\MCO_k$ prime to $6\MFm$, G.\,Robert proved that $\psi\left(1;\MFm,\MFa^{-1}\MFm\right) \in k(\MFm)$, where $k(\MFm)$ is the ray class field of $k$ modulo $\MFm$.
More precisely, $\psi\left(1;\MFm,\MFa^{-1}\MFm\right) \in \MCO_{k(\MFm)}^\times$ if $\MFm$ is divisible by at least two distinct primes, and if $\MFm=\MFr^n$ with $\MFr$ a prime ideal and $n\in\MBN^\ast$, then $\psi\left(1;\MFm,\MFa^{-1}\MFm\right)$ is a unit outside of the primes above $\MFr$.
For any maximal ideal $\MFq$ of $\MCO_k$, prime to $\MFa$, by \cite[Corollaire 1.3, (ii-1)]{robert89} we have
\begin{equation}
\label{ellipticunit}
N_{k(\MFm\MFq)/k(\MFm)} 
\left( \psi\left(1;\MFm\MFq,\MFa^{-1}\MFm\MFq\right) \right) ^{w_\MFm / w_{\MFm\MFq}} =
\lA
\begin{array}{lcl}
\psi\left(1;\MFm,\MFa^{-1}\MFm\right)^{1-\left(\MFq,k(\MFm)/k\right)^{-1}} & \mathrm{if} & \MFq\nmid\MFm,\\
&&\\
\psi\left(1;\MFm,\MFa^{-1}\MFm\right) & \mathrm{if} & \MFq\mid\MFm, 
\end{array}
\right.
\end{equation}
where $\left(\MFq,k(\MFm)/k\right)$ is the Fr\"obenius of $\MFq$ in $k(\MFm)/k$, and $w_\MFm$ is the number of roots of unity in $k$ which are congruent to $1$ modulo $\MFm$.
Moreover, by \cite[Corollaire 1.3, (v-1)]{robert89} we have 
\begin{equation}
\label{ellipticunitcongruence}
\psi\left(1;\MFm\MFq,\MFa^{-1}\MFm\MFq\right) \equiv \psi\left(1;\MFm,\MFa^{-1}\MFm\right)^{\left(\MFq,k(\MFm)/k\right)} \quad \text{modulo $(\MFq)_{\MFm\MFq}$},
\end{equation}
where $(\MFq)_{\MFm\MFq}$ is the product of the prime ideals in $\MCO_{k\left(\MFm\MFq\right)}$ above $\MFq$.

\begin{df}
Let $F\subseteq\MBC$ be a finite abelian extension of $k$, and write $\mu(F)$ for the group of roots of unity in $F$.
We write $\Psi_F$ for the $\MBZ\left[\Gal(F/k)\right]$-submodule of $F^\times$ generated by $\mu(F)$ and by all the norms
\[N_{k(\MFm)/k(\MFm)\cap F} \left( \psi\left(1;\MFm,\MFa^{-1}\MFm\right) \right),\]
where $\MFm$ is a nonzero proper ideal of $\MCO_k$ and $\MFa$ is any nonzero ideal of $\MCO_k$ prime to $6\MFm$.
Then, we define the group 
\[C_F \, := \, \Psi_F\cap\MCO_F^\times.\]
\end{df}

\begin{rmq}\label{leopoldt}
For any $n\in\MBN$, $\MCU_n$ is canonically identified to the pro-$p$-completion of the group of semi-local units $U_n$ of $K_n$.
Hence the natural inclusions $\MCO_{K_n}^\times \hookrightarrow U_n$ induce norm compatible canonical maps $\MCE_n\rightarrow\MCU_n$.
The Leopoldt conjecture, which is known to be true for abelian extensions of $k$, states that this map is injective.
Taking the projective limits, we obtain a natural injection $\MCE_\infty \hookrightarrow \MCU_\infty$.
\end{rmq}

\begin{pr}\label{torsionetGGGnfini}
The $\MBZ_p\left[\left[T\right]\right]$-module $\MCU_\infty / \MCC_\infty$ is finitely generated and torsion.
\end{pr}

\noindent\textsl{Proof.}
For all $n\in\MBN$, we let $St_n$ be the group of Stark units defined in \cite[Definition 3.2]{jilali-oukhaba}, and we set $\overline{St}_n := \MBZ_p \otimes_\MBZ St$ and $\overline{St}_\infty := \varprojlim\,\overline{St}_n$ (projective limit with respect to the norm maps).
It is well known that Stark units can be constructed by means of elliptic units (for instance, see \cite[Chapitre V,4]{oukhaba91} for a precise statement).
Then it is an easy matter to verify that $St_n\subseteq C_n$ for all $n\in\MBN$.
Hence $\overline{St}_\infty \subseteq \MCC_\infty$, and we just have to show that $\MCU_\infty / \overline{St}_\infty$ is finitely generated and torsion.
By \cite[Theorem 3.2 and Proposition 2.1]{jilali-oukhaba}, we know that $\overline{St}_\infty$ is torsion-free of rank $\left[ K_0 : k \right]$ over $\MBZ_p[[T]]$.
Then from \cite[Theorem 25]{iwasawa73} and Remark \ref{leopoldt}, we deduce that $\MCU_\infty / \overline{St}_\infty$ is finitely generated and torsion over $\MBZ_p[[T]]$.
\hfill $\square$
\\

\section{Euler systems.}\label{sectioneuler}

Let us write $A_k$ as a direct product of cyclic $p$-groups,
\[A_k \quad = \quad \la\cl\left(\MFp_1\right)\ra \times \cdots \times \la\cl\left(\MFp_r\right)\ra,\]
where $\MFp_1$, ..., $\MFp_r$ are prime ideals of $\MCO_k$, prime to $p$, and $\cl\left(\MFp_i\right)$ is the class of $\MFp_i$ in $\Cl\left(\MCO_k\right)$.
For any $i\in\{1,...,r\}$, let $p^{\SCr_i}$ be the order of $\la\cl(\MFp_i)\ra$, and we choose $\GGa_i\in\MCO_k$ be such that $\GGa_i\MCO_k = \MFp_i^{p^{\SCr_i}}$.
Let $\SCr:=\sum_{i=1}^r\SCr_i$ and let $\SCm\neq1$ be a power of $p$, such that $p^\SCr = \#(A_k) \leq \SCm$.
Let $\GGw:=1$ if $p\neq2$, and $\GGw:=-1$ if $p=2$.

We denote by $\MCL_F$ the set of maximal ideals $\ell$ of $\MCO_k$ such that $\ell$ splits completely in $F\left( \mu_\SCm, \sqrt [\SCm] {\GGw}, \sqrt[\SCm]{\GGa_1},...,\sqrt[\SCm]{\GGa_r}\right)/k$, and such that $\ell\notin\{\MFp_1,...,\MFp_r\}$.
We denote by $\MCS_F$ the set of squarefree ideals of $\MCO_k$ whose prime divisors belongs to $\MCL_F$.
As in \cite[Lemma 3.1]{oukhaba-viguie10}, we define for each $\ell\in\MCL_F$ a cyclic subextension $F\left(\ell\right)$ of $k(\ell)F$, of degree $\SCm$, which is totally ramified above $\ell$ and unramified anywhere else.
For $\MFm:=\ell_1\cdots\ell_n$ an ideal in $\MCS_F$, we define $F\left(\MFm\right) := F\left(\ell_1\right)\cdots F\left(\ell_n\right)$, the compositum of the fields $F\left(\ell_i\right)$.

For any ideal $\MFm\neq0$ of $\MCO_k$, we denote by $\MCS_F(\MFm)$ the set of ideals in $\MCS_F$ which are prime to $\MFm$.
We denote by $\MSU_F(\MFm)$ the set of maps $\GGe:\MCS_F(\MFm) \rightarrow (k^{ab})^\times$ satisfying the conditions (a) to (d) below.
\\

(a) $\GGe(\MFa)\in F\left(\MFa\right)^\times$ for all $\MFa\in\MCS_F(\MFm)$.
\\

(b) $\GGe(\MFa)\in\MCO_{F\left(\MFa\right)}^\times$ if $\MFa\neq(1)$.
\\

(c) $N_{F\left(\MFa\ell\right)/F\left(\MFa\right)}\left(\GGe(\MFa\ell)\right) = \GGe(\MFa)^{\left(\ell,F\left(\MFa\right)/k\right)-1}$ for all $\MFa\in\MCS_F(\MFm)$ and all $\ell\in\MCL_F$ which is prime to $\MFm\MFa$.
\\

(d) $\GGe(\MFa\ell) \equiv \GGe(\MFa)^{\left( N(\ell)-1 \right) / \SCm}$ modulo all prime ideals of $\MCO_{F\left(\MFa\right)}$ above $\ell$.
\\

\begin{rmq}\label{eulersystemstart}
Let $\MSU := \mathop{\cup} \MSU_F(\MFm)$, where the union is over all nonzero ideal $\MFm$ of $\MCO_k$.
Then for any $u\in C(F)$, there exists $\GGe\in\MSU$ such that $\GGe(1)=u$ (see \cite[Proposition 1.2]{rubin91}).
\end{rmq}

For any ideal $\MFa\neq(0)$ of $\MCO_k$ and any $\GGe\in\MSU_F(\MFa)$, we denote $\GGk_\GGe: \MCS_F(\MFa) \rightarrow F^\times/ \left(F^\times\right)^\SCm$ the map defined as in \cite[Proposition 2.2]{rubin91}.
For $\ell\in\MCL_F$, we let $\MCI_{F,\ell} := \ba{\oplus}{\GGl\mid\ell} \MBZ\GGl$ be the free $\MBZ$-module generated by the prime ideals of $\MCO_F$ lying above $\ell$.
For any $x\in F^\times$, we denote by $(x)_\ell\in\MCI_{F,\ell}$ and $[x]_\ell\in\MCI_{F,\ell}/\SCm\MCI_{F,\ell}$ the projections of the fractional ideal $(x):=x\MCO_K$.
Consider the map
\[\GGt_\ell : F\left(\ell\right)^\times \longrightarrow \left(\MCO_F/\ell\MCO_F\right)^\times / \left(\left(\MCO_F/\ell\MCO_F\right)^\times\right)^\SCm,\]
which associates to $z$ the sum $\oplus_{\GGl\mid\ell}z_\GGl$ such that the image of $z^{1-\GGs_\ell}$ in $\left(\MCO_F/\GGl\right)^\times$ is equal to $\left(z_\GGl\right)^{\left(N(\ell)-1\right) /\SCm}$.
As in \cite[Proposition 2.3]{rubin91}, there exists a unique $\Gal(F/k)$-equivariant isomorphism
\[\GVf_\ell: \left(\MCO_F/\ell\MCO_F\right)^\times / \left(\left(\MCO_F/\ell\MCO_F\right)^\times\right)^\SCm \rightarrow \MCI_{F,\ell} /\SCm \MCI_{F,\ell},\]
satisfying the relation $\left(\GVf_\ell\circ\GGt_\ell\right)(x) = \left[N_{ F\left(\ell\right) / F } (x)\right]_\ell$.
For $x\in F^\times$, we can choose $y\in F\left(\ell\right)^\times$ such that $xy^\SCm$ is a unit at the prime ideals of $\MCO_{F\left(\ell\right)}$ above $\ell$.
We denote by $\lA xy^\SCm\rA$ the class of $xy^\SCm$ in $\left(\MCO_{F\left(\ell\right)} / \ell'\right)^\times / \left(\left(\MCO_{F\left(\ell\right)} / \ell'\right)^\times\right)^\SCm  \simeq  \left(\MCO_F / \ell\MCO_F\right)^\times / \left(\left(\MCO_F / \ell\MCO_F\right)^\times\right)^\SCm$, where $\ell'$ is the product of the prime ideals of $\MCO_{F(\ell)}$ above $\ell$.
Then we set $\GVf_\ell(x) := \GVf_\ell\left(\lA xy^\SCm\rA\right)$, which does not depend on the choice of $y$.

Then as in \cite[Proposition 2.4]{rubin91}, for any ideal $\MFa\neq(0)$ of $\MCO_k$, any $\GGe\in\MSU_F(\MFa)$, any $\MFm\in\MCS_F(\MFa)$ with $\MFm\neq(1)$, and any maximal $\ell$ of $\MCO_k$, we have
\begin{equation}
\label{logelletphi}
\left[ \GGk_\GGe(\MFm) \right]_\ell = \lA
\begin{array}{lcl}
0 & \textrm{if} & \ell\nmid\MFm, \\
\GVf_\ell\left( \GGk_\GGe\left(\MFm\ell^{-1}\right) \right) & \textrm{if} & \ell\mid\MFm. \\
\end{array}
\right.
\end{equation}
For any $x\in F^\times$, we denote by $\la x\ra_\SCm$ the class of $x$ in $F^\times/\left(F^\times\right)^\SCm$.
For any $n\in\MBN$ we denote by $\mu_{n}$ the group of $n$-th roots of unity in $\MBC$.
We set $\mu_{p^\infty} := \mathop{\cup}_{n=0}^\infty\mu_{p^n}$.
For any extension $L\subseteq F$ of $k$ and any maximal ideal $\MFq$ of $\MCO_L$, we denote by $v_\MFq$ the normalized valuation at $\MFq$, and by $\bar{v}_\MFq: L^\times / \left(L^\times\right)^\SCm \rightarrow \MBZ / \SCm\MBZ$ the map defined from $v_\MFq$ by taking the quotient.

The following theorem is a classical step in the Euler system machinery.
The first versions are due to Rubin (see \cite[Theorem 3.1]{rubin91}), and to Greither for abelian extensions over $\MBQ$ (see \cite[Theorem 3.7]{greither92}).
We follow the proof of Bley (see \cite[Theorem 3.4]{bley06}), with slight modifications to cover the case $p\vert\#\left(\Cl\left(\MCO_k\right)\right)$.

\begin{theo}\label{theoeulergreither}
Let $\MFf$ be the conductor of $F/k$, and set $c:=v_{\bar{\MFp}}\left(\MFf\right)$.
We set $G_F:=\Gal(F/k)$.
Assume that we are given an ideal class $\MFc\in A_F$, a finite $\MBZ\left[G_F\right]$-submodule $W$ of $F^\times/ \left(F^\times\right)^\SCm$, and a $G_F$-morphism $\Psi : W \rightarrow \MBZ/\SCm\MBZ \left[G_F\right]$.

Assume that for all $w\in W$, all $i\in\{1,...,r\}$, and all prime $\MFq$ of $F$ above $\MFp_i$, $\bar{v}_{\MFq}(w) = 0$.
Assume also that for any $i=1,...,r$, $\MFp_i$ is unramified in $F/k$.
Let $m$ be a positive integer divisible by $p^{2c+1}$.
Then there are infinitely many maximal ideals $\GGl$ of $\MCO_F$ such that

\noindent (i) $\cl_p(\GGl) = \MFc^m$.

\noindent (ii) $\ell:=\GGl\cap\MCO_k$ belongs to $\MCL_F$.

\noindent (iii) For all $w\in W$, $[w]_\ell = 0$.

\noindent (iv) There exists $u\in\left(\MBZ/\SCm\MBZ\right)^\times$, such that for all $w\in W$, $\GVf_\ell(w) = up^{3c+\SCr+4} \Psi\left(w\right)\GGl$.
\end{theo}

\noindent\textsl{Proof.}
Let $H_F$ be the Hilbert $p$-class field of $F$.
Let
\[F_i:=\lA \begin{array}{lll} F\left(\mu_\SCm\right) & \text{if} & i=0, \\
F_{i-1}\left( \sqrt[\SCm]{\GGa_i} \right) & \text{if} & 1\leq i\leq r. \end{array} \right.\]
Exactly as in \cite[proof of Theorem 3.4]{bley06}, one can prove the following claims.

\noindent Claim (A) $\left[H_F\cap F\left(\mu_{p^\infty}\right) : F\right]\leq p^c$.

\noindent Claim (B) $\Gal\left(H_F\cap F_r\left( \sqrt[\SCm]{\GGw}, \sqrt[\SCm]{W} \right) / F\right) $ is annihilated by $p^{2c+1}$.

\noindent Claim (C) The cokernel of the canonical map from Kummer theory
\[ \MFK : \Gal\left(F_0 \left( \sqrt[\SCm]{W} \right) / F_0 \right) \hookrightarrow \Hom \left( W, \mu_\SCm \right)\]
 is annihilated by $p^{c+2}$.

Let us remark that $F_{i-1}\left(\sqrt[\SCm]{W}\right)/F_{i-1}$ is unramified at $\MFp_i$ since by hypothesis $\SCm\mid v_{\MFp_i}(w)$ for all $\la w\ra_\SCm \in W$.
On the other hand $[F_i:F_{i-1}]$ divides $\SCm$ and the ramification index of $\MFp_i$ in $F_i/F_{i-1}$ is at least $\SCm p^{-\SCr_i}$.
Therefore
\begin{equation}
p^{\SCr_i}\quad \text{annihilates} \quad \Gal\left( F_{i-1}\left(\sqrt[\SCm]{W}\right) \bigcap F_i / F_{i-1}\right).
\label{sdjfkdfs}
\end{equation}
Let $L_i:=F_0\left(\sqrt[\SCm]{W}\right)\cap F_i$.
As $L_i\cap F_{i-1}=L_{i-1}$ we have
\begin{equation}
\label{OYhHgF}
\Gal\left(L_i/L_{i-1}\right) \simeq \Gal\left(L_iF_{i-1}/F_{i-1}\right).
\end{equation}
Since $\Gal\left(L_iF_{i-1}/F_{i-1}\right)$ is a quotient of $\Gal\left( F_{i-1}\left(\sqrt[\SCm]{W}\right) \bigcap F_i / F_{i-1}\right)$, this implies that $p^{\SCr_i}$ annihilates $\Gal\left(L_i/L_{i-1}\right)$ thanks to (\ref{sdjfkdfs}).
In particular $p^\SCr$ annihilates $\Gal\left(L_r/F_0\right)$, and we deduce Claim (D) below.

\bigskip

\noindent Claim (D) $\Gal\left( F_0\left( \sqrt[\SCm]{W}\right) \bigcap F_r\left( \sqrt[\SCm]{\GGw} \right) / F_0\right)$ is annihilated by $p^{\SCr+1}$.

\bigskip

Let $\GGz$ be a primitive $\SCm$-root of unity, and $\GGi: \MBZ/\SCm\MBZ \left[G_F\right] \rightarrow \mu_\SCm$ be the group morphism such that $\GGi(\GGs)=0$ for $\GGs\in G_F\setminus\{1\}$ and $\GGi(1)=\GGz$.
Combining Claim (C) and Claim (D), one may find 
$\GGa\in \Gal\left( F_r\left( \sqrt[\SCm]{\GGw}, \sqrt[\SCm]{W} \right) / F_0 \right)$ such that
\begin{equation}
\label{MOJbhGF}
\GGa_{\vert F_r\left(\sqrt[\SCm]{\GGw}\right)}=1 \quad \text{and} \quad \MFK\left(\GGa_{\vert F_0\left(\sqrt[\SCm]{W} \right)}\right) = \left(\GGi\circ\Psi\right)^{p^{\SCr+c+3}}.
\end{equation}
From Claim (B), we may choose $\GGb\in\Gal \left( H_FF_r\left(\sqrt[\SCm]{\GGw}, \sqrt[\SCm]{W}\right) / F \right)$ such that
\begin{equation}
\label{MohhGGhff}
\GGb_{\vert F_r\left( \sqrt[\SCm]{\GGw}, \sqrt[\SCm]{W} \right)}=\GGa^{p^{2c+1}} \quad \text{and} \quad \GGb_{\vert H_F} = \MFc^m.
\end{equation}
Now, from (\ref{MOJbhGF}) we see that $\GGb\in\Gal\left(H_FF_r\left( \sqrt[\SCm]{\GGw}, \sqrt[\SCm]{W} \right) / F_r\left(\sqrt[\SCm]{\GGw}\right)\right)$.

By the \v Cebotarev density theorem, we can find infinitely many primes $\GGl$ in $\MCO_F$, of absolute degree $1$, prime to $\prod_{i=1}^r\MFp_i$, 
such that $\GGl\cap\MCO_k$ is unramified in $H_FF_r\left( \sqrt[\SCm]{\GGw}, \sqrt[\SCm]{W}\right)/k$, 
and such that the conjugacy class of $\GGb$ in $\Gal \left( H_FF_r\left(\sqrt[\SCm]{\GGw}, \sqrt[\SCm]{W}\right) / F \right)$ is the Fr\"obenius of $\GGl$.
Then condition (i) of Theorem \ref{theoeulergreither} holds as a consequence of the general properties of the Fr\"obenius.
The condition (ii) is also satisfied since $\GGb$ is the identity on $F_r\left(\sqrt[\SCm]{\GGw}\right)$.
Let $w\in W$.
Then for any prime $\GGl'$ of $\MCO_{F_0\left(\sqrt[\SCm]{W}\right)}$ above $\GGl$, we have $\bar{v}_\GGl\left(w\right) = \bar{v}_{\GGl'}\left(w\right) = \SCm \bar{v}_{\GGl'}\left(\sqrt[\SCm]{w}\right) = 0$, and condition (iii) follows.
Condition (iv) is proved as in the proof of \cite[Theorem 8.1,(iii)]{rubin91}.
\hfill $\square$
\\

For any $\MBZ_p[G]$-module $M$ and any $m\in M$, we denote by $m_\chi$ the canonical image of $m$ in $M_\chi$.

\begin{lm}\label{exceedingly}
Let $\MCG$ be a subgroup of $G_F$, and let $\chi$ be an irreducible $\MBC_p$ character of $\MCG$.
Let $\ell_1,...,\ell_i\in\MCL_F$, and for any $j=1,...,i$, let $\GGl_j$ be a prime of $\MCO_F$ above $\ell_j$, and let $\cl_p(\GGl_j)$ be the image of $\GGl_j$ in $A_F$.
Let $x\in F^\times$ be such that $v_\MFq(x)\in\SCm\MBZ$ for any prime $\MFq$ of $\MCO_F$ which is prime to $\ell_1\cdots\ell_i$.
Let $W$ be the $\MBZ_p\left[G_F\right]$-span of the image of $x$ in $F^\times / \left( F^\times \right)^\SCm$, and let $L$ be the $\MBZ_p \left[ G_F \right]$-module of $A_F$ generated by $\cl_p(\GGl_1),...,\cl_p(\GGl_{i-1})$.
Assume that there are $Z,g,\eta\in\MBZ_p\left[G_F\right]$ such that
\\

\noindent $\mathrm{(i)}$ $Z.\Ann_{\MBZ_p \left[ G_F \right]_\chi} \left( \left[ \cl_p(\GGl_i) \right]_{L,\chi} \right) \subseteq g\MBZ_p\left[ G_F \right]_\chi$, where $\Ann_{\MBZ_p \left[ G_F \right]_\chi} \left( \left[ \cl_p(\GGl_i) \right]_{L,\chi} \right)$ is the annihilator of the image $\left[ \cl_p(\GGl_i) \right]_{L,\chi}$ of $\cl_p(\GGl_i)$ in $\left(A_F/L\right)_\chi$.
\\

\noindent $\mathrm{(ii)}$ $\MBZ_p\left[ G_F \right]_\chi / g\MBZ_p\left[ G_F \right]_\chi$ is finite.
\\

\noindent $\mathrm{(iii)}$ $\#\left( \eta \left( \left( \MCI_{F,\ell_i} / \SCm \MCI_{F,\ell_i} \right) / W' \right)_\chi \right) \#\left(A_{F,\chi}\right) \leq \SCm$, where $W'$ is the image of $W$ in $\MCI_{F,\ell_i} / \SCm \MCI_{F,\ell_i}$ through $w\mapsto[w]_{\ell_i}$.
\\

Then, there exists a morphism of $\MBZ_p\left[ G_F \right]$-modules 
\[\Psi : W_\chi \rightarrow \MBZ/\SCm\MBZ \left[ G_F \right]_\chi\]
such that 
\[g\Psi\left(\la x\ra_{\SCm,\chi}\right)\GGl_{i,\chi} = Z\eta [x]_{\ell_i,\chi}.\]
\end{lm}

\noindent\textsl{Proof.}
We refer the reader to \cite[Lemma 3.12]{greither92}.
\hfill $\square$

\section{The ideal class group.}

\begin{pr}\label{Ainfty}
The projective limit $A_\infty$ is a finitely generated torsion $\MBZ_p[[T]]$-module.
Moreover, for all $n\in\MBN$, $A_{\infty,\GGG_n}$ and $A_\infty^{\GGG_n}$ are finite.
\end{pr}

\noindent\textsl{Proof.}
We refer the reader to \cite[Proof of Theorem 1.4]{rubin88}.
\hfill $\square$
\\

Let $\tau:A_{\infty,\chi} \rightarrow \op{\oplus}{j=1}{s} \GGL_\chi/P_j$ be a pseudo-isomorphism of $\GGL_\chi$-modules, where $P_1,...,P_s$ are nonzero polynomials in $\GGL_\chi$.

Let $M_\infty:=\varprojlim\left(M_n\right)$ be a $\MBZ_p[[T]]$-module, projective limit of $\MBZ_p\left[\GGG/\GGG_n\right]$-modules $M_n$.
For all $n\in\MBN$, we denote by $\Ker_nM_\infty$ and $\Cok_nM_\infty$ the respective kernel and cokernel of the canonical map $M_{\GGG_n}\rightarrow M_n$.

\begin{lm}\label{taun}
There is $c_3\in\MBN$, and for all $n\in\MBN$, there is a morphism of $\GGL_\chi$-modules
\[\tau_n:A_{n,\chi} \rightarrow \op{\oplus}{j=1}{s} \GGL_\chi/(P_j,1-\GGg_n)\]
such that $\Cok\left(\tau_n\right)$ is annihilated by $p^{c_3}$.
\end{lm}

\noindent\textsl{Proof.}
Let $m\in\MBN$ be such that $K_\infty/K_m$ is totally ramified above $\MFp$.
By \cite[Lemma 13.15]{washington97}, there is a $\MBZ_p[[\GGG_m]]$-submodule $Y$ of $A_\infty$ such that for all $n\geq m$,
the canonical map $A_\infty \rightarrow A_n$ induces an isomorphism
\begin{equation*}
A_\infty/ \nu_{m,n}Y \; \xymatrix{\ar[r]^-{\sim} &} \; A_n,
\end{equation*}
where $\nu_{m,n}\in\MBZ_p[[T]]$ is defined by $\nu_{m,n} := \left(1-\GGg_n\right) / \left(1-\GGg_m\right)$.
Therefore for all $n\geq m$, we have $\Cok_nA_\infty=0$ and 
\begin{equation}
\label{MnoiZA}
\Ker_nA_\infty \; \simeq \; \nu_{m,n}Y \, / \, \left(1-\GGg_n\right)A_\infty.
\end{equation} 
Multiplication by $\nu_{m,n}$ induces a surjection
\begin{equation}
\label{NvfYg}
Y/\left(1-\GGg_m\right)A_\infty \; \xymatrix{\ar@{->>}[r] &} \nu_{m,n}Y \, / \, \left(1-\GGg_n\right)A_\infty,
\end{equation}
from which we deduce that for all $n\geq m$, $\Ker_nA_\infty$ is a quotient of $\Ker_mA_\infty$.
Since $\Ker_nA_\infty$ is finite, by Proposition \ref{Ainfty} we see that the orders of $\Ker_nA_\infty$ and $\Cok_nA_\infty$ are bounded independantly of $n$.
We choose $\GGa\in\MBN$ such that for all $n\in\MBN$, $p^\GGa$ annihilates $\Ker_nA_\infty$ and $\Cok_nA_\infty$.
On the other hand, since $\Cok_nA_\infty=0$, we have the exact sequence below for any $n\geq m$,
\begin{equation*}
\xymatrix{
\left(\Ker_nA_\infty\right)_\chi \ar[r] & \left(A_\infty\right)_{\GGG_n,\chi} \ar[r] & A_{n,\chi} \ar[r] & 0.
}\end{equation*}
This shows that $p^\GGa$ annihilates $\Ker\left(\left(A_{\infty,\chi}\right)_{\GGG_n} \rightarrow A_{n,\chi}\right)$, for all $n\geq m$.
Moreover for all $n\in\MBN$, $\left(A_{\infty,\chi}\right)_{\GGG_n} \simeq \left(A_{\infty,\GGG_n}\right)_\chi$ is finite from Proposition \ref{Ainfty}.
Thus we may choose $\GGa$ such that $p^\GGa$ also annihilates $\Ker\left(\left(A_{\infty,\chi}\right)_{\GGG_n} \rightarrow A_{n,\chi}\right)$ for all $n\in\MBN$.
Then choose $\GGb\in\MBN$ such that $p^\GGb$ annihilates $\Cok(\tau)$, and set $c_3:=2\GGa+\GGb$.
Let $\bar{\tau}_n: \left(A_{\infty,\chi}\right)_{\GGG_n} \rightarrow \op{\oplus}{j=1}{s}\GGL_\chi/(P_j,1-\GGg_n)$ be the morphism of $\GGL_\chi$-modules defined from $\tau$ by taking the quotients, and set
\[\tau_n:A_{n,\chi} \xymatrix{\ar[r] &} \op{\oplus}{j=1}{s}\GGL_\chi/(P_j,1-\GGg_n), \; x \xymatrix{\ar@{|->}[r] &} p^\GGa\bar{\tau}_n(y),\]
where $y\in\left(A_{\infty,\chi}\right)_{\GGG_n}$ is such that its image in $A_{n,\chi}$ is $p^\GGa x$.
It is straightforward that $\tau_n$ is well-defined, and that the condition of the lemma is satisfied.
\hfill $\square$

\section{Global units.}\label{sectionprel}

Let us fix $\chi$ an irreducible $\MBC_p$ character of $G$, and let $\psi : G\rightarrow\MBZ_p$ the irreducible $\MBQ_p$ character of $G$ such that $\chi\vert\psi$.
Also we denote by $\MFu_\chi$ a fixed uniformizer of $\MBZ_p(\chi)$.

\begin{lm}\label{fchi}
There is a finite set $I$, a family $(n_i)_{i\in I}\in\MBN^I$, and a pseudo-isomorphism of $\GGL_\chi$-modules:
\begin{equation}
\label{equafchi}
\GGT: \MCE_{\infty,\chi} \rightarrow \GGL_\chi\oplus\ba{\oplus}{i\in I}\left(\GGL_\chi/\MFu_\chi^{n_i}\right).
\end{equation}
\end{lm}

\noindent\textsl{Proof.}
From \cite[Theorem 25]{iwasawa73} and Remark \ref{leopoldt}, we know that $\MCE_\infty$ is finitely generated over $\MBZ_p[[T]]$.
We denote by $\GGL(\psi)$ the principal ring $e_\psi \left( \MBQ_p \otimes_{\MBZ_p} \GGL \right)$.
From the tautological exact sequence $0 \rightarrow \MCE_\infty \rightarrow \MCU_\infty \rightarrow \MCU_\infty / \MCE_\infty \rightarrow 0$ we deduce the following exact sequence of $\GGL(\psi)$-modules
\begin{equation}
\label{IGJHffj}
0 \rightarrow e_\psi\left( \MBQ_p \otimes_{\MBZ_p} \MCE_\infty \right) \rightarrow e_\psi\left( \MBQ_p \otimes_{\MBZ_p} \MCU_\infty \right) \rightarrow e_\psi\left( \MBQ_p \otimes_{\MBZ_p} \MCU_\infty / \MCE_\infty \right) \rightarrow 0.
\end{equation}
But we know that $e_\psi\left( \MBQ_p \otimes_{\MBZ_p} \MCU_\infty \right)$ is free of rank $1$ over $\GGL(\psi)$, thanks to Proposition \ref{QptensorZpsemilocalunits}.
Moreover, $e_\psi\left( \MBQ_p \otimes_{\MBZ_p} \MCU_\infty / \MCE_\infty \right)$ is $\GGL(\psi)$-torsion by Proposition \ref{torsionetGGGnfini}.
As $\GGL(\psi)$ is a principal ring, we see from (\ref{IGJHffj}) that $e_\psi\left( \MBQ_p \otimes_{\MBZ_p} \MCE_\infty \right)$ is free of rank $1$ over $\GGL(\psi)$.
The isomorphisms
\[e_\psi\left( \MBQ_p \otimes_{\MBZ_p} \MCE_\infty \right) \simeq \MBQ_p \otimes_{\MBZ_p} \MCE_{\infty,\chi} \quad\text{and}\quad \GGL(\psi) \simeq \MBQ_p \otimes_{\MBZ_p} \GGL_\chi\]
imply that the $\GGL_\chi$-torsion of $\MCE_{\infty,\chi}$ is annihilated by some power of $p$, and the $\GGL_\chi$-rank of $\MCE_{\infty,\chi}$ is $1$.
\hfill $\square$
\\

\begin{lm}\label{lemmaczeroexists}
There is $\left(c_0,n_0\right)\in\MBN^2$ such that for all $n\in\MBN$, $p^{c_0}(\GGg_{n_0}-1)$ annihilates $\Cok_n\MCE_\infty$ and $\Ker_n\MCE_\infty$.
\end{lm}

\noindent\textsl{Proof.}
The proof is similar to \cite[Corollary 3.9]{oukhaba10}.
\hfill $\square$
\\

For every $n\in\MBN$, the projection $\MCE_\infty\rightarrow\MCE_n$ induces a natural map $\pi_{n,\chi} : \left(\MCE_{\infty,\chi}\right)_{\GGG_n} \rightarrow \MCE_{n,\chi}$. 

\begin{lm}\label{lemmaczero}
For all $n\in\MBN$, $p^{2c_0}(\GGg_{n_0}-1)^2$ annihilates $\Cok\left(\pi_{n,\chi}\right)$ and $\Ker\left(\pi_{n,\chi}\right)$.
\end{lm}

\noindent\textsl{Proof.} 
Let $n\in\MBN$ and let $\MCT:=\Tor_{\MBZ_p[G]}^1\left(\Cok_n\MCE_\infty,\MBZ_p(\chi)\right)$.
We denote by $\widetilde{\MCE}_n$ the image of $\MCE_\infty$ in $\MCE_n$, and we write $\widetilde{\pi}_{n,\chi}: \left(\MCE_{\infty,\chi}\right)_{\GGG_n} \rightarrow \widetilde{\MCE}_{n,\chi}$.
From the following commutative exact diagram,
\[\xymatrix{
 & & \MCT \ar[d] & \\
\left(\Ker_n\MCE_\infty\right)_\chi \ar[d] \ar[r] & \left(\MCE_{\infty,\chi}\right)_{\GGG_n} \ar[d]^{\pi_{n,\chi}} \ar[r] & \widetilde{\MCE}_{n,\chi} \ar[d] \ar[r] & 0 \\
0 \ar[r] & \MCE_{n,\chi} \ar@{=}[r] & \MCE_{n,\chi} \ar[d] \ar[r] & 0 \\
 & & \left(\Cok_n\MCE_\infty\right)_\chi \ar[d] & \\
 & & 0 & ,\\
}\]
we deduce an exact sequence
\begin{equation}
\label{nbJTBjy}
\xymatrix{\left(\Ker_n\MCE_\infty\right)_\chi \ar[r] & \Ker\left(\pi_{n,\chi}\right) \ar[r] & \widetilde{\MCT} \ar[r] & 0, }
\end{equation}
where $\widetilde{\MCT}$ is the image of $\MCT$ in $\widetilde{\MCE}_{n,\chi}$.
By Lemma \ref{lemmaczeroexists}, we know that $p^{c_0}(\GGg_{n_0}-1)$ annihilates $\MCT$ and $\left(\Ker_n\MCE_\infty\right)_\chi$.
Therefore (\ref{nbJTBjy}) implies the desired result for $\Ker\left(\pi_{n,\chi}\right)$.
On the other hand, since $\Cok\left(\pi_{n,\chi}\right)\simeq\left(\Cok_n\MCE_\infty\right)_\chi$ the lemma is entirely proved.
\hfill $\square$
\\

Let $\mathrm{pr}:\GGL_\chi\oplus\ba{\oplus}{i\in I}\left(\GGL_\chi/\MFu_\chi^{n_i}\right) \rightarrow \GGL_\chi$ be the canonical projection, and for every $n$ let
\[\GGT_n : \left(\MCE_{\infty,\chi}\right)_{\GGG_n} \xymatrix{\ar[r] &} \MBZ_p(\chi)\left[\GGG/\GGG_n\right]\]
be the map obtained from $\PR\circ\GGT$ by taking the quotients.
By Remark \ref{leopoldt}, $\MCE_\infty / \MCC_\infty$ is a submodule of $\MCU_\infty / \MCC_\infty$, hence is finitely generated and torsion over $\MBZ_p[[T]]$ by Proposition \ref{torsionetGGGnfini}.
We denote by $h_\chi$ a generator of $\Char\left(\MCE_\infty/\MCC_\infty\right)_\chi$.
The following lemma is the analogue of \cite[Lemma 3.5]{bley06}.

\begin{lm}\label{GVtnchi}
Let $n\in\MBN$.
Then there is a map
\[\GVt_{n,\chi} : \MCE_{n,\chi} \rightarrow \MBZ_p(\chi)\left[\GGG/\GGG_n\right], x \mapsto \left(\GGg_{n_0}-1\right)^2 p^{2c_0} \GGT_n \left(\underline{x}\right),\]
where $\underline{x}\in \MCE_{\infty,\chi}$ is such that $\left(\GGg_{n_0}-1\right)^2 p^{2c_0} x = \pi_{n,\chi} \left(\underline{x}\right)$.
Then there are $\left(\nu,c_1,c_2\right)\in\MBN^3$ and $h'_\chi\in\GGL_\chi$ such that

$\mathrm{(i)}$ $h'_\chi\vert h_\chi$ in $\GGL_\chi$.

$\mathrm{(ii)}$ For all $n\in\MBN$, $h'_\chi$ is prime to $1-\GGg_n$ in $\GGL_\chi$.

$\mathrm{(iii)}$ For all $n\in\MBN$, $\left(\GGg_{\nu}-1\right)^{c_1}p^{c_2} h'_\chi \MBZ_p(\chi)\left[\GGG/\GGG_n\right] \subseteq \GVt_{n,\chi} \left(\IM \left(\MCC_{n,\chi}\right) \right)$,
where $\IM \left(\MCC_{n,\chi}\right)$ is the image of $\MCC_{n,\chi}$ in $\MCE_{n,\chi}$.
\end{lm}

\noindent\textsl{Proof.}
One may use Lemma \ref{lemmaczero} to verify that $\GVt_{n,\chi}$ is well-defined.
We leave the details to the reader.
The module $h_\chi\cdot\left(\MCE_\infty/\MCC_\infty\right)_\chi$ is finite, and so is $h_\chi\cdot\left(\GGT\left(\MCE_{\infty,\chi}\right) / \GGT\left(\IM\left(\MCC_{\infty,\chi}\right)\right)\right)$.
Since $\Cok\left(\mathrm{pr}\circ \GGT\right)$ is also finite, we can choose $m\in\MBN$ such that $p^mh_\chi \in \mathrm{pr}\circ \GGT\left(\IM\left(\MCC_{\infty,\chi}\right)\right)$.
Let $z\in\IM\left(\MCC_{\infty,\chi}\right)$ be such that $p^mh_\chi = \mathrm{pr}\circ \GGT(z)$.
Then in $\MBZ_p(\chi)\left[\GGG/\GGG_n\right]$ we have
\begin{equation}
p^{m+4c_0}\left(\GGg_{n_0}-1\right)^4h_\chi = p^{4c_0}\left(\GGg_{n_0}-1\right)^4 \GGT_n(z) = \GVt_{n,\chi}\left(\pi_{n,\chi}(z)\right). \label{ggyeqo}
\end{equation}
Let $\MCQ$ be the set of prime ideals $\MFq$ of $\GGL_\chi$ of height $1$, and for any $\MFq\in\MCQ$, let $P_\MFq$ be a generator of $\MFq$.
Since $\GGL_\chi$ is factorial, there is a unit $u\in\GGL_\chi^\times$ and a family $(n_\MFq)_{\MFq\in\MCQ}\in\MBN^\MCQ$ with finite support such that $h_\chi=u\prod_{\MFq\in\MCQ}P_\MFq^{n_\MFq}$.
We set $h'_\chi := \prod_{\MFq\in\MCQ'}P_\MFq^{n_\MFq}$, where $\MCQ'$ is the set of all $\MFq\in\MCQ$ such that $\MFq$ is prime to $1-\GGg_n$, for all $n\in\MBN$.
Since $1-\GGg_n$ divides $1-\GGg_{n+1}$ for all $n\in\MBN$, we can choose $\nu\in\MBN$ and $c_1\in\MBN$ such that $\left(\GGg_{n_0}-1\right)^4h_\chi$ divides $\left(\GGg_{\nu}-1\right)^{c_1}h'_\chi$.
Also, we set $c_2:=m+4c_0$.
Then the lemma follows from (\ref{ggyeqo}).
\hfill $\square$
\\

\section{Proof of Theorem \ref{mainconj}.}\label{sectiomainconj}

This section is devoted to the proof of Theorem \ref{mainconj}.
We choose $c_0$ and $n_0$ as in Lemma \ref{lemmaczeroexists}, $c_1\geq2$, $c_2$ and $\nu$ as in Lemma \ref{GVtnchi}, and $c_3$ as in Lemma \ref{taun}.
Let us define 
\[d := 3v_{\bar{\MFp}} (\MFf) + \SCr + 4 \quad\text{and}\quad \GGD_i := p^{(i-2) (c_3+2d) +d+c_2} [K_0:k]^{i-1} \quad(i\geq2),\]
where $\MFf$ is a nonzero ideal of $\MCO_k$ such that $K_\infty\subseteq\mathop{\cup}_{n=0}^\infty k\left(\MFf\MFp^n\right)$.
Let $n\in\MBN$.
Since $h'_\chi$ is prime to $1-\GGg_n$, the factor group 
\[\MBZ_p[G_n]_\chi/\GGD_{s+1} h'_\chi\MBZ_p[G_n]_\chi \; \simeq \; \GGL_\chi/\left( (1-\GGg_n)\GGL_\chi + \GGD_{s+1} h'_\chi \GGL_\chi \right)\] 
is finite.
Let $\SCm$ be a power of $p$ such that
\begin{equation}
\label{choixdeSCm}
\#A_k \, \#A_{n,\chi} \, \#\left(\MBZ_p[G_n]_\chi / \GGD_{s+1}h'_\chi\MBZ_p[G_n]_\chi\right) \; \leq \; \SCm.
\end{equation}
Let us also introduce the following notation.
If $\GGl$ is a maximal ideal of $\MCO_{K_n}$ such that $\ell:=\GGl\cap\MCO_k \in\MCL_{K_n}$, then we denote by $\GGw_\GGl$ and $\bar{\GGw}_\GGl$ the maps
\[\GGw_\GGl: K_n^\times \xymatrix{\ar[r] &} \MBZ_p\left[G_n\right], \quad\text{such that}\quad \GGw_\GGl(x)\GGl=(x)_\ell,\]
and
\[\bar{\GGw}_\GGl: K_n^\times / \left(K_n^\times\right)^\SCm \xymatrix{\ar[r] &} \left(\MBZ/\SCm\MBZ\right)\left[G_n\right], \quad\text{such that}\quad \bar{\GGw}_\GGl\left(\la x\ra_\SCm\right)\GGl=[x]_\ell.\]
We know by Lemma \ref{taun} that for every $j\in\{1,...,s\}$, there is a class $\MFc_j\in A_n$ such that 
\[\tau_n\left(\MFc_{j,\chi}\right) = \left(0, ..., 0, p^{c_3}, 0, ..., 0\right),\] 
where $p^{c_3}$ is at the $j$-th place.
We recall that $\MFc_{j,\chi}$ is the image of $\MFc_j$ in $A_{n,\chi}$.
We also choose arbitrarily one more class $\MFc_{s+1}\in A_n$.
By the above Lemma \ref{GVtnchi} (iii), there is $\xi\in C_n$ such that  
\begin{equation}
\label{dudfngle}
\GVt_{n,\chi}(\xi') = \left(\GGg_{\nu}-1\right)^{c_1} p^{c_2} h'_\chi \quad\text{in}\quad \left(\MBZ/\SCm\MBZ[G_n]\right)_\chi,
\end{equation}
where $\xi'$ is the image of $\xi$ in $\IM\left(\MCC_{n,\chi}\right)$.
Let us now fix an ideal $\MFm$ of $\MCO_k$ and $\GVe\in\MSU_{K_n}(\MFm)$ such that $\GGk_\GVe(1)=\xi$.
This is possible thanks to Lemma \ref{eulersystemstart}.
The main step is to define recursively maximal ideals $\GGl_1,...,\GGl_{s+1}$ of $\MCO_{K_n}$ and ideals $\MFa_1,...,\MFa_{s+1}$ of $\MCO_k$ such that 
\\ 

\noindent(a) $\ell_i:=\GGl_i\cap\MCO_k$ belongs to $\MCL_{K_n}$ for all $i=1,...,s+1$.
\\

\noindent(b) $\cl_p\left(\GGl_i\right)=\MFc_i^{p^d}$ for all $i=1,...,s+1$.
\\

\noindent(c) $\MFa_i:=\ell_1\cdots\ell_i$.
\\

\noindent(d) $\bar{\GGw}_{\GGl_1} \left(\GGk_\GVe(\ell_1)\right)_\chi = u_1p^{c_2+d} [K_0:k]\left(\GGg_{\nu}-1\right)^{c_1}h'_\chi$ in $\left(\MBZ/\SCm\MBZ[G_n]\right)_\chi$, for some $u_1\in\left(\MBZ/\SCm\MBZ\right)^\times$.
\\

\noindent(e) For every $i\in\{2,...,s+1\}$ there is $u_i\in\left(\MBZ/\SCm\MBZ\right)^\times$ such that
\[P_{i-1}\bar{\GGw}_{\GGl_i}\left(\GGk_\GVe(\MFa_i)\right) = u_ip^{c_3+2d} [K_0:k]\left(\GGg_{\nu}-1\right)^{c_1^{i-1}} \bar{\GGw}_{\GGl_{i-1}} \left(\GGk_\GVe(\MFa_{i-1})\right).\]

Since we will use Theorem \ref{theoeulergreither} for $F:=K_n$, we assume from now until the end of this paper that
\[\text{the prime ideals $\MFp_1,...,\MFp_r$ are unramified in $K_\infty/k$.}\]
Let us consider the map $\GVp\circ\GVt_{n,\chi}\circ\eta:\MCO_{K_n}^\times\rightarrow\MBZ_p\left[G_n\right]$, 
where $\eta:\MCO_{K_n}^\times\rightarrow\MCE_{n,\chi}$ is the natural map and $\GVp:\MBZ_p(\chi)\left[\GGG/\GGG_n\right] \rightarrow \MBZ_p\left[G_n\right]$ 
is defined by $\GVp\left(\chi(g)\GGu\right):=\left[K_0:k\right]e_\psi g\GGu$ for all $g\in G$ and $\GGu\in\GGG/\GGG_n$.
Further, by taking the quotients we obtain a map
\[\Psi_1: \MCO_{K_n}^\times / \left(\MCO_{K_n}^\times\right)^\SCm \xymatrix{\ar[r] &} \MBZ/\SCm\MBZ\left[G_n\right].\]
Let $W_1$ be the $\MBZ_p\left[G_n\right]$-span of $\la\xi\ra_\SCm$.
We apply Theorem \ref{theoeulergreither} to the data
\[m:=p^d, \quad W:=W_1, \quad \Psi:=\Psi_1, \quad\text{and}\quad \MFc:=\MFc_1,\]
We obtain a maximal ideal $\GGl_1$ of $\MCO_{K_n}$ and $u_1\in\left(\MBZ/\SCm\MBZ\right)^\times$ such that $\cl_p(\GGl_1) = \MFc_1^{p^d}$, such that the ideal $\ell_1:=\GGl_1\cap\MCO_k$ belongs to $\MCL_{K_n}$, and such that for all $w\in W_1$, we have $[w]_{\ell_1}=0$ and 
\begin{equation}
\GVf_{\ell_1}(w) = u_1p^d \Psi_1\left(w\right)\GGl_1.
\end{equation}
We denote by $\bar{\GVt}_{n,\chi} : \MCE_{n,\chi} \rightarrow \MBZ/\SCm\MBZ[G_n]_\chi$ the morphism obtained from $\GVt_{n,\chi}$ by taking the quotients.
Then from (\ref{logelletphi}), we have
\begin{equation}
\label{sdidjfnfyg}
\left[\GGk_\GVe(\ell_1)\right]_{\ell_1} = \GVf_{\ell_1}(\xi) =  u_1p^d \Psi_1\left(\la\xi\ra_\SCm\right)\GGl_1 = u_1p^d \GVp \circ \bar{\GVt}_{n,\chi} \left( \xi' \right)\GGl_1
\end{equation}
in $\MCI_{K_n,\ell_1}/\SCm\MCI_{K_n,\ell_1}$.
From (\ref{sdidjfnfyg}) and (\ref{dudfngle}), we deduce that in $\left(\MBZ/\SCm\MBZ[G_n]\right)_\chi$ we have
\begin{eqnarray}
\bar{\GGw}_{\GGl_1} \left(\GGk_\GVe(\ell_1)\right) 
& = & u_1p^d [K_0:k]\bar{\GVt}_{n,\chi}\left(\xi'\right) \nonumber \\
& = & u_1p^{c_2+d} [K_0:k]\left(\GGg_{\nu}-1\right)^{c_1}h'_\chi.
\label{xidnfyhha}
\end{eqnarray}
Let $i\in\{2,...,s+1\}$, and assume that $\GGl_1,...,\GGl_{i-1}$ has been constructed.
From (d) and (e) we deduce
\begin{equation}
\label{equaalacon}
\left(\prod_{j=1}^{i-2}P_j\right) \bar{\GGw}_{\GGl_{i-1}} \left(\GGk_\GVe(\MFa_{i-1})\right) = \left(\prod_{j=1}^{i-1}u_j\right) \GGD_i \left(\GGg_{\nu}-1\right)^{c_1+\sum_{j=1}^{i-2}c_1^j} h'_\chi
\end{equation}
in $\left(\MBZ/\SCm\MBZ\left[G_n\right]\right)_\chi$, with the convention that an empty product is $1$ and an empty sum is $0$.

\begin{lm}\label{onpeututiliserexceedingly}
Let $L_i$ be the $\MBZ_p[G_n]$-submodule of $A_n$ generated by $\cl_p(\GGl_1),...,\cl_p(\GGl_{i-2})$, and let $W_i$ be the $\MBZ_p[G_n]$-span of the image of $\GGk_\GVe(\MFa_{i-1})$ in $K_n^\times / \left( K_n^\times \right)^\SCm$.
We set $\eta_i:=\left(\GGg_{\nu}-1\right)^{c_1^{i-1}}$, $Z_i:=p^{d+c_3}$, and we choose $g_i\in\MBZ_p[G_n]$ such that the image of $g_i$ and the image of $P_{i-1}$ in $\MBZ_p[G_n]_\chi$ are the same.
Then
\\

\noindent $\mathrm{(i)}$ $v_\MFq\left( \GGk_\GVe(\MFa_{i-1}) \right) \in \SCm\MBZ$ for all maximal ideal $\MFq$ of $\MCO_{K_n}$ which is prime to $\MFa_{i-1}$.
\\

\noindent $\mathrm{(ii)}$ $Z_i.\Ann_{\MBZ_p \left[ G_n \right]_\chi} \left( \left[ \cl_p(\GGl_{i-1}) \right]_{L_i,\chi} \right) \subseteq g_i\MBZ_p\left[ G_n \right]_\chi$, where $\left[ \cl_p(\GGl_{i-1}) \right]_{L_i,\chi}$ is the image of $\cl_p(\GGl_{i-1})$ in $\left(A_n/L_i\right)_\chi$.
\\

\noindent $\mathrm{(iii)}$ $\MBZ_p\left[ G_n \right]_\chi / g_i\MBZ_p\left[ G_n \right]_\chi$ is finite.
\\

\noindent $\mathrm{(iv)}$ $\#\left( \eta_i \left( \left( \MCI_{K_n,\ell_{i-1}} / \SCm \MCI_{K_n,\ell_{i-1}}\right) / W'_i \right)_\chi \right) \#\left(A_{n,\chi}\right) \leq \SCm$, where $W'_i$ is the image of $W_i$ in $\MCI_{K_n,\ell_{i-1}} / \SCm \MCI_{K_n,\ell_{i-1}}$ through $w\mapsto[w]_{\ell_{i-1}}$.
\end{lm}

\noindent\textsl{Proof.}
$\mathrm{(i)}$ is a direct consequence of (\ref{logelletphi}).
We have $\left(A_{\infty,\chi}\right)_{\GGG_n} \simeq \left(A_\infty\right)_{\GGG_n,\chi}$, and $\left(A_\infty\right)_{\GGG_n,\chi}$ is finite by Proposition \ref{Ainfty}.
Hence $g_i$ is prime to $1-\GGg_n$, whence  $\mathrm{(iii)}$.

Let $\GGa\in\Ann_{\MBZ_p \left[ G_n \right]_\chi} \left( \left[ \cl_p(\GGl_{i-1}) \right]_{L_i,\chi} \right)$.
We can define from $\tau_n$ a morphism of $\MBZ_p[G_n]_\chi$-modules
\[\tau'_n : \left(A_n/L_i\right)_\chi \rightarrow \MBZ_p[G_n]_\chi / g_i\MBZ_p\left[ G_n \right]_\chi,\]
such that the diagram below commutes
\[\xymatrix{
A_{n,\chi} \ar[rr]^-{\tau_n} \ar@{->>}[d] & & \op{\oplus}{j=1}{s}\GGL_\chi/(P_j,1-\GGg_n) \ar@{->>}[d]^{\phi} & \\
\left(A_n/L_i\right)_\chi \ar[rr]^-{\tau'_n} & & \MBZ_p[G_n]_\chi/g_i\MBZ_p[G_n]_\chi &,
}\]
where $\phi$ is the canonical projection
\[\op{\oplus}{j=1}{s}\GGL_\chi/(P_j,1-\GGg_n) \longrightarrow \GGL_\chi/\left(P_{i-1},1-\GGg_n\right) \simeq \MBZ_p[G_n]_\chi/g_i\MBZ_p[G_n]_\chi.\]
Then $\tau'_n\left(\MFc_{i-1,\chi}\right)^{p^d\GGa} = 0$, i.e $p^{d+c_3}\GGa\in g_i\MBZ_p\left[ G_n \right]_\chi$, so $\mathrm{(ii)}$ is verified.
From (\ref{equaalacon}), and since $c_1+\sum_{j=1}^{i-2}c_1^j\leq c_1^{i-1}$ (because $2\leq c_1$), we see that 
$\eta_i \left( \left(\MCI_{K_n,\ell_{i-1}} / \SCm \MCI_{K_n,\ell_{i-1}}\right) / W'_i \right)_\chi$ is cyclic over $\MBZ/\SCm\MBZ[G_n]_\chi$, annihilated by $\GGD_ih'_\chi$.
The condition (\ref{choixdeSCm}) then implies $\mathrm{(iv)}$.
\hfill $\square$
\\

Let us apply Lemma \ref{exceedingly} to the material furnished in Lemma \ref{onpeututiliserexceedingly}.
There is a morphism of $\MBZ_p[G_n]$-modules $\Psi'_i : W_{i,\chi} \rightarrow \MBZ/\SCm\MBZ[G_n]_\chi$ such that
\begin{equation}
\label{ehydgc}
g_i\Psi'_i \left(\la\GGk_\GVe\left(\MFa_{i-1}\right)\ra_{\SCm,\chi}\right) \GGl_{i-1,\chi} = Z_i\eta_i \left[ \GGk_\GVe\left(\MFa_{i-1}\right) \right]_{\ell_{i-1},\chi}.
\end{equation}

We define $\Psi_i$ by composing $\GVp\circ\Psi'_i$ with $W_i\rightarrow W_{i,\chi}$.
From Lemma \ref{onpeututiliserexceedingly}, $(i)$, we can apply Theorem \ref{theoeulergreither} to the data
\[F:=K_n, \quad m:= p^d, \quad \MFc:=\MFc_i, \quad W:=W_i, \quad\text{and}\quad \Psi:=\Psi_i.\]
There are a maximal ideal $\GGl_i$ of $\MCO_{K_n}$ and $u_i\in\left(\MBZ/\SCm\MBZ\right)^\times$ such that $\cl_p(\GGl_i) = \MFc_i^{p^d}$ (condition (b)), such that $\ell_i:=\GGl_i\cap\MCO_k$ belongs to $\MCL_{K_n}$ (condition (a)), and such that for all $w\in W_i$, $[w]_{\ell_i} = 0$ and $\GVf_{\ell_i}(w) = u_ip^d \Psi_i\left(w\right)\GGl_i$.
By (\ref{logelletphi}) we have
\begin{equation}
\label{equationyg}
\left[ \GGk_\GVe \left(\MFa_i\right) \right]_{\ell_i,\chi} = \GVf_{\ell_i} \left( \GGk_\GVe \left(\MFa_{i-1}\right) \right)_\chi = u_ip^d \Psi_i\left(\la\GGk_\GVe \left(\MFa_{i-1}\right)\ra_\SCm\right)\GGl_{i,\chi}
\end{equation}
in $\left(\MCI_{K_n,\ell_i} /\SCm \MCI_{K_n,\ell_i}\right)_\chi$.
Then in $\MBZ/\SCm\MBZ[G_n]_\chi$, by (\ref{equationyg}) and (\ref{ehydgc}) we have
\[\begin{array}{ccl}
P_{i-1} \bar{\GGw}_{\GGl_i} \left( \GGk_\GVe\left(\MFa_i\right) \right) & = & u_i [K_0:k] p^dg_i\Psi'_i\left(\la \GGk_\GVe\left(\MFa_{i-1}\right) \ra_{\SCm,\chi}\right) \\
 & = & u_i [K_0:k] p^dZ_i\eta_i \bar{\GGw}_{\GGl_{i-1}} \left( \GGk_\GVe\left(\MFa_{i-1}\right) \right), \\
\end{array}\]
which demonstrates (e).
So we can construct recursively the primes $\GGl_1,...,\GGl_{s+1}$, and from (d) and (e) we deduce
\begin{equation}
\label{equaalacondeux}
\left(\prod_{j=1}^{s}P_j\right) \bar{\GGw}_{\GGl_{s+1}} \left(\GGk_\GVe(\MFa_{s+1})\right) = \left(\prod_{j=1}^{s+1}u_j\right) \GGD_{s+2} \left(\GGg_{\nu}-1\right)^{c_1+\sum_{j=1}^s c_1^j} h'_\chi
\end{equation}
in $\MBZ/\SCm\MBZ[G_n]_\chi$.

By letting $n$ and $\SCm$ vary, this implies that $\prod_{j=1}^{s}P_j$ divides $\GGD_{s+2} \left(\GGg_{\nu}-1\right)^{c_1+\sum_{j=1}^s c_1^j} h'_\chi$ in $\GGL_\chi$.
By Proposition \ref{Ainfty} and since $\left(A_{\infty,\chi}\right)_{\GGG_{\nu}} \simeq \left(A_{\infty,\GGG_{\nu}}\right)_\chi$, $\Char\left(A_{\infty,\chi}\right)$ is prime to $(\GGg_{\nu}-1)$.
Then we deduce 
\begin{equation}
\label{DiViSiBiLiTy}
\Char\left(A_{\infty,\chi}\right) | \GGD_{s+2} \Char\left(\MCE_\infty/\MCC_\infty\right)_\chi,
\end{equation} 
which proves the assertion (ii) of Theorem \ref{mainconj}.
Now assume $p\notin\{2,3\}$.
Recall that class field theory gives an exact sequence
\begin{equation}
\label{tUiLEBLUE}
\xymatrix{
0 \ar[r] & \MCE_\infty/\MCC_\infty \ar[r] & \MCU_\infty/\MCC_\infty \ar[r] & B_\infty \ar[r] & A_\infty \ar[r] & 0,
}
\end{equation}
where $B_\infty$ is the Galois group of $\GGW_\infty/K_\infty$, whith $\GGW_\infty$ the maximal abelian pro-$p$-extension of $K_\infty$ which is unramified outside of the primes above $\MFp$.
Moreover by a result of Gillard (\cite[3.4. Th\'eor\`eme]{gillard85}), the $\mu$-invariant of $B_\infty$ over $\MBZ_p[[T]]$ vanishes,
\begin{equation}
\label{nulgillard}
\mu\left(B_\infty\right) \quad = \quad 0.
\end{equation}
Let us denote by $A_{\MFf,\infty}$, $\MCE_{\MFf,\infty}$, $\MCC_{\MFf,\infty}$, ..., the various objects attached to $K_{\MFf,\infty} := \mathop{\cup}_{n=0}^\infty k\left(\MFf\MFp^n\right)$.
By \cite[2.1 Theorem, p.\,109]{deshalit87}, the divisibility (\ref{DiViSiBiLiTy}) and (\ref{nulgillard}) applied to $K_{\MFf,\infty}$ implies that
\begin{equation}
\label{lOybf}
\Char\left(\left(A_{\MFf,\infty}\right)_\xi\right) = \Char\left(\MCE_{\MFf,\infty}/\MCC_{\MFf,\infty}\right)_\xi \quad\text{and}\quad \mu_\xi\left(\MCU_{\MFf,\infty}/\MCC_{\MFf,\infty}\right)_\xi = 0,
\end{equation}
for all irreducible $\MBC_p$-character $\xi$ of the torsion subgroup $G_\MFf$ of $\Gal\left( K_{\MFf,\infty} / k \right)$, and where $\mu_\xi$ is the $\mu$-invariant over $\MBZ_p(\xi)$.
Since $\left(\MCU_\infty/\MCC_\infty\right)_\chi$ is a quotient of $\left(\MCU_{\MFf,\infty}/\MCC_{\MFf,\infty}\right)_{\tilde{\chi}}$, 
where $\tilde{\chi}$ is the character of $G_\MFf$ defined by $\chi$, we deduce from (\ref{lOybf}) that $\mu_\chi\left(\MCU_\infty/\MCC_\infty\right)_\chi = 0$.
By (\ref{nulgillard}), the exact sequence (\ref{tUiLEBLUE}) gives
\begin{equation}
\label{lumPf}
\mu_\chi\left(\MCE_\infty/\MCC_\infty\right)_\chi = \mu_\chi\left(\MCU_\infty/\MCC_\infty\right)_\chi = 0 = \mu_\chi\left(A_{\infty,\chi}\right).
\end{equation}
By decomposing $H:=\Gal\left(K_{\MFf,\infty}/K_\infty\right)$ into a direct product of cyclic subgroups, we are reduced to the case where $H$ itself is cyclic.
Then classical arguments (see \cite[section 5]{rubin91}) show that $\left(A_{\MFf,\infty}\right)_H$ is pseudo-isomorphic to $A_\infty$, and we deduce
\begin{equation}
\label{HuyVJYt}
\GGl\left(\left(A_{\MFf,\infty}\right)_H\right) = \GGl\left(A_\infty\right).
\end{equation}
The cokernel of the norm map $\left(\MCE_{\MFf,\infty}/\MCC_{\MFf,\infty}\right)_H \rightarrow \MCE_\infty/\MCC_\infty$ is annihilated by $\#(H)$, hence $\GGl\left(\MCE_\infty/\MCC_\infty\right)\leq\GGl\left(\MCE_{\MFf,\infty}/\MCC_{\MFf,\infty}\right)_H$.
Together with (\ref{lOybf}) and (\ref{HuyVJYt}), it implies that
\begin{equation}
\label{HumalabarYt}
\GGl\left(\MCE_\infty/\MCC_\infty\right) \leq \GGl\left(\MCE_{\MFf,\infty}/\MCC_{\MFf,\infty}\right)_H = \GGl\left(\left(A_{\MFf,\infty}\right)_H\right) = \GGl\left(A_\infty\right).
\end{equation}
Finally the assertion (i) of Theorem \ref{mainconj} follows from (\ref{HumalabarYt}), (\ref{lumPf}) and (\ref{DiViSiBiLiTy}).

We draw the attention of the reader to our papers \cite{viguie11a} and \cite{viguie11c}, where we prove that a raw version of Theorem \ref{mainconj} $\mathrm{(i)}$ holds also for $p\in\{2,3\}$.

\bibliographystyle{amsplain}

\end{document}